\documentclass[12pt]{article}
\usepackage{latexsym,amsfonts,amssymb}
\usepackage[dvips]{color}
\usepackage{colordvi,multicol}
\usepackage{amsmath}
\usepackage{graphicx}
\usepackage{pdflscape}
\usepackage{rotating}

\usepackage{latexsym,amsfonts,amssymb}
\usepackage[dvips]{color}
\usepackage{colordvi,multicol}
\usepackage{amsmath}
\usepackage{graphicx}
\usepackage{pdflscape}
\usepackage{rotating}
\usepackage{breqn}

\usepackage{scalerel}
\usepackage{mathtools, amssymb}

\DeclareMathOperator{\p}{\mathbb{P}}

\DeclareMathOperator{\var}{\text{Var}}
\DeclareMathOperator{\prd}{\mathcal{P}(\mathbb{R}^d)}

\makeatletter
\newcommand{\mylabel}[2]{#2\def\@currentlabel{#2}\label{#1}}
\makeatother

\def \cal{\mathcal}

\newtheorem{thm}{Theorem}[section]
\newtheorem{cor}[thm]{Corollary}
\newtheorem{lem}[thm]{Lemma}
\newtheorem{pro}[thm]{Proposition}
\newtheorem{defn}[thm]{Definition}
\newtheorem{rem}[thm]{Remark}
\newtheorem{exa}[thm]{Example}

\date{}
\begin{document}

\title{\bf Existence of Periodic and Stationary Solutions to Distribution-Dependent SDEs}
 \author{Wei Sun and Ethan Wong\\ \\ \\
  {\small Department of Mathematics and Statistics}\\
    {\small Concordia University, Canada}\\ \\
{\small wei.sun@concordia.ca,\ \ \ \ ethan.wong@concordia.ca}}

\maketitle

\begin{abstract}

\noindent We investigate the periodic and stationary solutions of distribution-dependent stochastic differential equations. While generally, the semigroups associated with the equations are nonlinear, we show that the methods of weak convergence and Lyapunov functions can be combined to give efficient criteria for the existence of periodic and stationary solutions. Concrete examples are presented to illustrate the novel criteria.
\end{abstract}


\noindent  {\it MSC2020 subject classifications:} 60H10; 37A50.

\section{Introduction and Main Results}

A distribution-dependent stochastic differential equation (DDSDE), also known as a McKean-Vlasov SDE, is an SDE whose coefficients depend on the distribution of the solution. It describes a stochastic system with a large number of interacting particles, whose evolution is determined by both the microcosmic site and the macroscopic distribution of the particles. The study of DDSDEs has a long history that can be traced back to Kac's work \cite{Kac} on chaos propagation and McKean's seminal paper \cite{McKean2} on large stochastic interacting particle systems. Since then, the equations have been widely investigated in the literature; see, e.g., \cite{Dawson, Gartner, Huang1, Huang2, LL1, LL2, LL3, Sznitman} for detailed developments on the subject. To date, the theory of DDSDEs has found applications across a broad range of disciplines, including, but not limited to, finance, game theory, communication networks, neuroscience, population dynamics, and economics (cf. the monograph \cite{CARMONA}).

The concept of stationary solutions plays a fundamental role in studying the long-time behavior of random dynamical systems. In the past decades, many works have been devoted to studying stationary solutions of DDSDEs. Here, we list some of them that are closely related to our paper. Ahmed and Ding \cite{AhmedDing}, Veretennikov \cite{V}, and Butkovsky \cite{B} gave different sufficient conditions for the existence and uniqueness of invariant measures and discussed convergence to those invariant measures for DDSDEs with additive noise. Wang \cite{wang2017distributiondependent} developed conditions for the existence and uniqueness of strong solutions to DDSDEs and studied the ergodicity of the solutions, with applications to specific models, including the homogeneous Landau equation. Bogachev et al. \cite{Boga} studied convergence in variation of probability solutions of nonlinear Fokker-Planck-Kolmogorov equations to stationary solutions and obtained sufficient conditions for the exponential convergence of solutions to the stationary solution. Zhang \cite{Zhang.nonuniqueness} and Bao et al. \cite{Bao2022} used fixed point theorems to show the existence of invariant measures for DDSDEs and functional DDSDEs, respectively. Additionally, Zhang considered the non-uniqueness of invariant measures with applications to the phase transition phenomenon in mathematical physics. By considering the linearization of nonlinear Fokker-Planck-Kolmogorov equations, Ren et al. \cite{ren2020linearization} characterized the ergodicity of solutions to DDSDEs. Liu and Ma \cite{LIU2022138.MarkovSwitching} imposed Lyapunov conditions to show the existence and uniqueness of invariant measures for DDSDEs with Markovian switching. Du et al. \cite{EmpiricalApprox} showed that, under a monotonicity condition, weighted empirical measures can approximate the invariant measures of DDSDEs.

In this paper, we investigate the existence of periodic and stationary solutions to DDSDEs. It is well known that periodic solutions have played an important role in the study of dynamical systems since the pioneering work of Poincar\'e \cite{Pon1, Pon2, Pon3, Pon4}. However, most of the existing literature on the long-time behavior of DDSDEs focuses on the time-homogeneous case. In the time-homogeneous case, a typical method to establish the existence of stationary solutions is to first freeze the distribution component of DDSDEs to take advantage of the ergodic theory for linear Markov semigroups, and then use the Schauder fixed point theorem (cf. \cite{Bao2022,Zhang.nonuniqueness}). It seems difficult, if not impossible, to apply this standard method to establish the existence of periodic solutions to time-inhomogeneous DDSDEs. Another method to study stationary solutions of DDSDEs is to couple the stochastic process with its distribution and then use the projection technique (cf. \cite{Liu2,ren2020linearization}). A drawback of this method is that very strong conditions must be imposed to obtain stationary solutions to the coupled SDEs, which are not necessary if we are mainly interested in the stationary solutions of DDSDEs. Our paper is among the first to study periodic solutions of time-inhomogeneous DDSDEs. The method developed here is novel and can also be applied to obtain new results for stationary solutions to time-homogeneous DDSDEs. Regarding periodic solutions of ordinary SDEs, we refer the reader to the monograph \cite{khasminskii2011stochastic} and \cite{Chen,Guo2019PeriodicSO,Guo2021HybridJump,Xu,XuHuang,ZhangX} for some previous works related to this paper. We also refer the reader to the recent preprint \cite{Bao2} which investigated random periodic solutions for DDSDEs. We would like to point out that the notion of random periodic solution considered in \cite{Bao2} is entirely different from the notion of periodic solution considered in our paper (cf. \cite[Remark 2.3]{Bao2}), and the main results given in \cite[Theorems 2.1 and 2.4]{Bao2} are very distinct from ours.

Let $\{B_t\}_{t\geq0}$ be a standard $d$-dimensional Brownian motion on a complete probability space with natural filtration $(\Omega, \mathcal{F}, \{\mathcal{F}_t\}_{t \geq 0}, \p)$. Denote by $\prd$ the space of all probability measures on $\mathbb{R}^d$ equipped with the weak topology, $\mathcal{L}_\eta$  the distribution of a random variable $\eta$ on $\mathbb{R}^d$, $|\cdot|$ and $\langle\cdot,\cdot\rangle$ the standard Euclidean norm and inner product, respectively,  and  $\|A\|$ the trace norm of a matrix $A$. Let \[
    b: [0,\infty)\times\mathbb{R}^d\times\prd \rightarrow \mathbb{R}^d,\quad  \quad \sigma: [0,\infty)\times\mathbb{R}^d\times\prd \rightarrow \mathbb{R}^d\otimes\mathbb{R}^d
\]
be measurable maps. We are interested in the DDSDE on $\mathbb{R}^d$:
\begin{equation}\label{DDSDE}
dX_t=b(t,X_t,{\cal L}_{X_t})dt+\sigma(t,X_t,{\cal L}_{X_t})dB_t.
\end{equation}
We restrict ourselves to the following subspace of $\prd$ for some $\vartheta \in [1,\infty)$:
\[
    \mathcal{P}_{\vartheta}(\mathbb{R}^d) := \left\{\mu \in \prd : \Vert\mu\Vert_{\vartheta} :=  \left(\int_{\mathbb{R}^d}|x|^\vartheta\mu(dx)\right)^{\frac{1}{\vartheta}} < \infty\right\},
\]
which is a Polish space under the Wasserstein-$\vartheta$ metric:
\[
    W_\vartheta(\mu,\nu) := \inf_{\pi \in \mathcal{C}(\mu,\nu)}\left(\int_{\mathbb{R}^d \times \mathbb{R}^d}|x-y|^\vartheta\pi(dx,dy)\right)^{\frac{1}{\vartheta}}, \quad\quad \mu, \nu \in \mathcal{P}_{\vartheta}(\mathbb{R}^d),
\]
where $\mathcal{C}(\mu,\nu)$ is the set of all couplings for $\mu$ and $\nu$.

If the DDSDE (\ref{DDSDE}) has strong existence and uniqueness in ${\mathcal{P}_{\vartheta}(\mathbb{R}^d)}$, then the solution $\{X_t\}_{t\ge0}$ is Markovian. For $s \geq 0$, let $\{X_{s,t}^\mu\}_{t \geq s}$ solve (\ref{DDSDE}) with initial distribution $\mathcal{L}_{X_{s,s}^\mu} = \mu\in {\mathcal{P}_{\vartheta}(\mathbb{R}^d)}$. By the uniqueness of solutions, we get
\begin{equation}\label{memoryless}
    X_{s,t}^\mu = X_{r,t}^{X_{s,r}^\mu}, \quad 0 \leq s \leq r \leq t.
\end{equation}
Further, if the DDSDE (\ref{DDSDE}) has also ${\mathcal{P}_{\vartheta}(\mathbb{R}^d)}$-weak  uniqueness, then we may define a semigroup
$\{P_{s,t}^*\}_{t\ge s}$ on ${\mathcal{P}_{\vartheta}(\mathbb{R}^d)}$ by letting $P_{s,t}^*\mu = {\cal L}_{X_{s,t}}$ for ${\cal L}_{X_{s,s}}=\mu\in {\mathcal{P}_{\vartheta}(\mathbb{R}^d)}$. In fact, by (\ref{memoryless}), we have that
$$
P_{s,t}^*=P_{r,t}^*P_{s,r}^*,\quad 0 \leq s \leq r \leq t.
$$
For simplicity of notation, denote $P_{t}^*\mu:=P_{0,t}^*\mu $ for $\mu\in {\mathcal{P}_{\vartheta}(\mathbb{R}^d)}$ and $t\ge0$.

A key obstacle to DDSDEs comes from the fact that the semigroup
$\{P_{s,t}^*\}_{t\ge s}$ is generally nonlinear. For example, let us consider the following simple DDSDE:
\begin{equation}\label{Dec5a}
    dX_t=\var(X_t)dt + dB_t.
\end{equation}
Suppose that $X_0$ is independent of $\{B_t\}_{t\geq0}$. Denote by $\delta_x$ the Dirac measure at point $x\in \mathbb{R}^d$ and $\nu_0$ the standard normal distribution on $\mathbb{R}^d$. Let $X_{t}^\mu:=X_{0,t}^\mu$ for $\mu\in\prd$ and $t\ge0$. Writing (\ref{Dec5a}) in the integral form:
\[
    X_t = X_0 + \int_0^t \var(X_s)ds + B_t,
\]
we have that
\begin{equation}\label{dec10b}
    X_t^{\delta_x} = x + \int_0^t \var(X_s^{\delta_x})ds + B_t, \quad\quad  X_t^{\nu_0} = X_0^{\nu_0} + \int_0^t \var(X_s^{\nu_0})ds + B_t.
\end{equation}
Taking variances on both sides, we get $\var(X_t^{\delta_x}) = t$ and $\var(X_t^{\nu_0}) = 1 + t$. Then, substituting back into (\ref{dec10b}), we obtain the following solutions to (\ref{Dec5a}):
$$
    X^{\delta_x}_t = x + \frac{t^2}{2} + B_t,\quad\quad
    X^{\nu_0}_t = X^{\nu_0}_0 + t +\frac{t^2}{2} + B_t.
$$
Thus,
\begin{equation}\label{Dec5c}
P^*_{t}\delta_x \sim {\rm N}\left(x + \frac{t^2}{2}, t\right),\quad\quad
P^*_{t}\nu_0 \sim {\rm N}\left( t + \frac{t^2}{2}, 1 + t\right).
\end{equation}

By (\ref{Dec5c}), we get
\begin{eqnarray*}
        \int_{-\infty}^\infty(P^*_t\delta_x)\nu_0(dx)
        &=& \int_{-\infty}^\infty\left(\frac{1}{\sqrt{2\pi t}}e^{-{\frac{(y - x - \frac{t^2}{2})^2}{2t} }}dy\right)\cdot\frac{1}{\sqrt{2\pi}} e^{-\frac{x^2}{2}}dx \\
        &=& \frac{1}{2\pi \sqrt{t}} \int_{-\infty}^\infty e^{-\frac{(1+t)x^2 - 2(y-\frac{t^2}{2}) x+ (y-\frac{t^2}{2})^2}{2t}} dydx  \\
               &=& \frac{1}{2\pi \sqrt{t}} \int_{-\infty}^\infty e^{-\frac{(\sqrt{1+t}\cdot x - \frac{y - \frac{t^2}{2}}{\sqrt{1+t}} )^2}{2t} } \cdot e^{- \frac{(y- \frac{t^2}{2})^2}{2(1+t)}} dxdy\\
               &=& \frac{1}{\sqrt{2\pi(1+t)}} e^{- \frac{(y- \frac{t^2}{2})^2}{2(1+t)}}dy.
\end{eqnarray*}
Hence,
$$
    \int_{-\infty}^\infty(P^*_{t}\delta_x)\nu_0(dx) \sim {\rm N}\left(\frac{t^2}{2}, 1+t\right)\not\sim P^*_{t}\nu_0.
$$
Therefore, the semigroup $\{P^*_{s,t}\}_{t\ge0}$ associated with the DDSDE (\ref{Dec5a}) is nonlinear.

To overcome difficulties caused by nonlinearity and time-inhomogeneouity of the semigroup $\{P^*_{s,t}\}_{t\ge0}$, we will combine the methods of weak convergence (cf. \cite{KurtzMarkov}) and Lyapunov functions (cf. \cite{khasminskii2011stochastic})  to give efficient criteria for the existence of periodic and stationary solutions to DDSDEs. Throughout this paper, we fix a $\theta > 0$.
\begin{defn} The stochastic process $\{X_t\}_{t\ge 0}$ is said to be a $\theta$-periodic solution to the DDSDE (\ref{DDSDE}) if  the following conditions hold.
    \begin{enumerate}
        \item $\{X_t\}_{t\ge 0}$ is a solution to the DDSDE (\ref{DDSDE}).
        \item For any $k \in \mathbb{N}$ and any sequence $0 \leq t_1 < t_2< \cdots< t_n$, the joint distribution of $X_{t_1 +k\theta}, X_{t_2 +k\theta}, \dots, X_{t_n +k\theta}$ is independent of $k$.
    \end{enumerate}
    \end{defn}

We put the following assumptions:
\begin{itemize}
    \item [(\textbf{H0})] The coefficients $b$ and $\sigma$ are $\theta$-periodic in $t$, i.e.,
    $$
    b(t,x,\mu)=b(t+\theta,x,\mu),\ \ \sigma(t,x,\mu)=\sigma(t+\theta,x,\mu),\ \
    \forall t\ge0, x\in\mathbb{R}^d, \mu\in{\cal P}(\mathbb{R}^d).
    $$
    \item[(\textbf{H1})] The DDSDE (\ref{DDSDE}) has strong existence and both strong and weak uniqueness in ${\mathcal{P}_{\vartheta}(\mathbb{R}^d)}$.
    \item[(\textbf{H2})] $P_{\theta}^*({\cal P}_{\vartheta}(\mathbb{R}^d))\subset {\cal P}_{\vartheta}(\mathbb{R}^d)$ and $P_{\theta}^*$ is a continuous map from the Polish space
        $({\cal P}_{\vartheta}(\mathbb{R}^d),W_{\vartheta})$ to itself.
\end{itemize}

\begin{rem} In \cite[Theorem 2.1]{wang2017distributiondependent}, Wang showed that  conditions (H1) and (H2) are implied by some assumptions on the continuity, monotonicity and growth of coefficients. Also, Wang presented the result \cite[Theorem 6.3]{wang2017distributiondependent} for the equivalence of the weak existence/uniqueness and the strong existence/uniqueness under an additional assumption.
\end{rem}

Let $C^{1,2}([0,\infty)\times\mathbb{R}^d)$ be the space of all real-valued functions on $[0,\infty)\times\mathbb{R}^d$ which are continuously differentiable in the first component, and twice continuously differentiable in the second component. For $V\in C^{1,2}([0,\infty)\times\mathbb{R}^d)$, $t\geq0$, $x\in\mathbb{R}^d$  and $\mu\in{\cal P}(\mathbb{R}^d)$, we define
\begin{eqnarray}\label{Lya Generator}
    (LV)(t,x,\mu)&:=&\partial_t V(t,x)+\langle b(t,x,\mu),\partial_x V(t,x)\rangle\nonumber\\
&&+\frac{1}{2}{\rm trace}\left(\sigma\sigma^T(t,x,\mu)\cdot\partial_x^2V(t,x)\right).
\end{eqnarray}

\begin{thm}\label{thm1}
    Let $\vartheta \in [1,\infty)$. Suppose that conditions (H0)--(H2) hold. Then, there exists a $\theta$-periodic solution to the DDSDE (\ref{DDSDE}) if one of the following three conditions holds.
\begin{itemize}
        \item [(\textbf{H3a})] There exist $V\in C^{1,2}([0,\infty)\times\mathbb{R}^d)$, $C_0,C_1,C_2>0$ and $r>\vartheta$ such that
    \[
        V(t,x)\ge -C_0,\quad \forall t\ge0, x\in\mathbb{R}^d,
    \]
    and
    \[
        (LV)(t,x,\mu)\le -C_1|x|^{r}+C_2(1+\|\mu\|^{\vartheta}_{\vartheta}),\quad \forall t\ge0, x\in\mathbb{R}^d, \mu\in{\cal P}_{\vartheta}(\mathbb{R}^d).
    \]
Moreover, $\int_0^t \|{\cal L}_{X^{\delta_0}_s}\|^{\vartheta}_{\vartheta}ds<\infty$ for any $t>0$.
\item [(\textbf{H3b})]
There exist $V\in C^{1,2}([0,\infty)\times\mathbb{R}^d)$, $C_0>0$, $C_1>C_2>0$, $C_3>0$  and $r>\vartheta$ such that
    \[
        V(t,x)\ge -C_0,\quad \forall t\ge0, x\in\mathbb{R}^d,
    \]
    and
    \[
        (LV)(t,x,\mu)\le -C_1|x|^r+C_2\|\mu\|^r_r+C_3,\quad \forall t\ge0, x\in\mathbb{R}^d, \mu\in{\cal P}_{r}(\mathbb{R}^d).
    \]
Moreover, $\int_0^t \|{\cal L}_{X^{\delta_0}_s}\|^r_rds<\infty$ for any $t>0$.
\item [(\textbf{H3c})] There exist $V\in C^{1,2}([0,\infty)\times\mathbb{R}^d)$, $C_1,C_2>0$ and $r>\vartheta$ such that
    \begin{equation*}\label{Oct14a}
    V(t,x)\ge V(0,x)\ge C_1|x|^{r} -C_2,\ \ \ \ \forall t\ge0, x\in\mathbb{R}^d,
    \end{equation*}
and
$$(LV)(t,x,\mu)\le 0,\ \ \ \ \forall t\ge0, x\in\mathbb{R}^d, \mu\in{\cal P}(\mathbb{R}^d).$$
\end{itemize}
\end{thm}

We may also consider the time-independent version of the DDSDE (\ref{DDSDE}):
\begin{equation}\label{timeIndDDSDE}
dX_t=b(X_t,{\cal L}_{X_t})dt+\sigma(X_t,{\cal L}_{X_t})dB_t,
\end{equation}
where \[
    b: \mathbb{R}^d\times\prd \rightarrow \mathbb{R}^d,\quad  \quad \sigma: \mathbb{R}^d\times\prd \rightarrow \mathbb{R}^d\otimes\mathbb{R}^d
\]
are measurable maps.

\begin{defn} The stochastic process $\{X_t\}_{t\ge 0}$ is said to be a stationary solution to the DDSDE (\ref{timeIndDDSDE}) if the following conditions hold.
    \begin{enumerate}
        \item $\{X_t\}_{t\ge 0}$ is a solution to the DDSDE (\ref{timeIndDDSDE}).
        \item For any $h>0$ and any sequence $0 \leq t_1 < t_2< \cdots< t_n$, the joint distribution of  $X_{t_1 +h}, X_{t_2 +h}, \dots, X_{t_n +h}$ is independent of $h$.
    \end{enumerate}
   \end{defn}

We put the following assumptions:
\begin{itemize}
    \item[(\textbf{H1'})] The DDSDE  (\ref{timeIndDDSDE})  has strong existence and both strong and weak  uniqueness in ${\mathcal{P}_{\vartheta}(\mathbb{R}^d)}$.
   \item [(\textbf{H2'})]  For any $t>0$, $P_t^*({\cal P}_{\vartheta}(\mathbb{R}^d))\subset {\cal P}_{\vartheta}(\mathbb{R}^d)$ and $P_t^*$ is a continuous map from the Polish space $({\cal P}_{\vartheta}(\mathbb{R}^d),W_{\vartheta})$ to itself.
\end{itemize}

Let $C^{2}(\mathbb{R}^d)$ be the space of all real-valued twice continuously differentiable functions on $\mathbb{R}^d$. For $V\in C^{2}(\mathbb{R}^d)$, $x\in\mathbb{R}^d$  and $\mu\in{\cal P}(\mathbb{R}^d)$, we define
\[
    (LV)(x,\mu):=\langle b(x,\mu),\partial_x V(x)\rangle+\frac{1}{2}{\rm trace}\left(\sigma\sigma^T(x,\mu)\cdot\partial_x^2V(x)\right).
\]

\begin{cor}\label{cor}
    Let $\vartheta \in [1,\infty)$. Suppose that conditions (H1') and (H2') hold. Then, there exists a stationary solution to the DDSDE (\ref{timeIndDDSDE}) if one of the following three conditions holds.
\begin{itemize}
        \item [(\textbf{H3a'})] There exist $V\in C^{2}(\mathbb{R}^d)$, $C_0,C_1,C_2>0$ and $r>\vartheta$ such that
    \[
        V(x)\ge -C_0,\quad \forall x\in\mathbb{R}^d,
    \]
    and
    \[
        (LV)(x,\mu)\le -C_1|x|^{r}+C_2(1+\|\mu\|^{\vartheta}_{\vartheta}),\quad \forall x\in\mathbb{R}^d, \mu\in{\cal P}_{\vartheta}(\mathbb{R}^d).
    \]
Moreover,  $\int_0^t \|{\cal L}_{X^{\delta_0}_s}\|^{\vartheta}_{\vartheta}ds<\infty$ for any  $t>0$.
\item [(\textbf{H3b'})] There exist $V\in C^{2}(\mathbb{R}^d)$, $C_0>0$,  $C_1>C_2>0$, $C_3>0$ and $r>\vartheta$ such that
    \[
        V(x)\ge -C_0,\quad \forall x\in\mathbb{R}^d,
    \]
    and
    \[
        (LV)(x,\mu)\le -C_1|x|^{r}+C_2\|\mu\|^{r}_{r}+C_3,\quad \forall x\in\mathbb{R}^d, \mu\in{\cal P}_{r}(\mathbb{R}^d).
    \]
Moreover,  $\int_0^t \|{\cal L}_{X^{\delta_0}_s}\|^r_rds<\infty$ for any  $t>0$.
\item [(\textbf{H3c'})] There exist $V\in C^{2}(\mathbb{R}^d)$, $C_1,C_2>0$ and $r>\vartheta$ such that
     \begin{equation*}
    V(x)\ge C_1|x|^{r} -C_2,\ \ \ \ \forall x\in\mathbb{R}^d,
    \end{equation*}
and
$$(LV)(x,\mu)\le 0,\ \ \ \ \forall x\in\mathbb{R}^d, \mu\in{\cal P}(\mathbb{R}^d).$$
\end{itemize}
\end{cor}

It is well known that a McKean-Vlasov SDE may have several stationary solutions (cf. \cite{Dawson,Duong,Tugaut,Zhang.nonuniqueness}). Different from many existing results, e.g., \cite[Theorem 3.1] {wang2017distributiondependent}, \cite[Theorem 4.1]{LIU2022138.MarkovSwitching} and \cite[Theorems 2.1 and 2.4]{Bao2}, our Theorem \ref{thm1} and Corollary \ref{cor} do not exclude the possibility of non-uniqueness of stationary or periodic solutions to DDSDEs. It is interesting to further establish the uniqueness and non-uniqueness criteria for periodic and stationary solutions to DDSDEs in the framework of Theorem \ref{thm1} and Corollary \ref{cor}. We will leave this for future work.

The rest of the paper is organized as follows. In section 2, we give the proof of Theorem \ref{thm1}. The proof of Corollary \ref{cor} is very similar, which we will omit. In section 3, we present concrete examples to demonstrate the main results.

\section{Proof of Theorem \ref{thm1}}\label{sec2}\setcounter{equation}{0}

In this section, we give the proof of Theorem \ref{thm1}.
\begin{defn} A probability measure $\mu$  on $\mathbb{R}^d$ is said to be a $\theta$-periodic measure of the DDSDE (\ref{DDSDE}) if $P_{\theta}^*\mu = \mu$.
    \end{defn}

We first present a lemma that describes the  relationship between a $\theta$-periodic solution and a $\theta$-periodic measure.

\begin{lem}\label{lemma}
 Suppose that conditions \textbf{(H0)} and \textbf{(H1)} hold and ${\cal L}_{X_0}\in {\mathcal{P}_{\vartheta}(\mathbb{R}^d)}$ for some $\vartheta \in [1,\infty)$. Then, $\{X_t\}_{t\ge0}$ is a $\theta$-periodic solution to the DDSDE (\ref{DDSDE})  if and only if ${\cal L}_{X_0}$ is a $\theta$-periodic measure of the DDSDE (\ref{DDSDE}).
\end{lem}

\noindent {\it Proof.} The proof of the necessity is obvious. We only prove the sufficiency.
Suppose that $\{X_t\}_{t\ge0}$ is the unique strong  solution to the DDSDE (\ref{DDSDE}) with ${\cal L}_{X_0}\in {\mathcal{P}_{\vartheta}(\mathbb{R}^d)}$ being a $\theta$-periodic measure. Set $\mu={\cal L}_{X_0}$.
    Then,  $\mathcal{L}_{X_{0}^\mu} = \mathcal{L}_{X_{\theta}^\mu} = \mu$.

    By performing the substitution $s = u+ \theta$, letting $\Bar{B}_t = B_{t+\theta} - \Bar{B}_\theta$, and using periodicity of the coefficients $b$ and $\sigma$, we get
    \begin{equation*}
        \begin{split}
            X_{t+\theta} &= X_\theta + \int_\theta^{t+\theta}b(s,X_s, \mathcal{L}_{X_s})ds + \int_\theta^{t+\theta}\sigma(s,X_s, \mathcal{L}_{X_s})dB_s\\
            &= X_\theta + \int_0^{t}b(u+\theta,X_{u+\theta}, \mathcal{L}_{X_{u+\theta}})du + \int_0^{t}\sigma(u+\theta,X_{u+\theta}, \mathcal{L}_{X_{u+\theta}})d\Bar{B}_u\\
            &= X_\theta + \int_0^{t}b(u,X_{u+\theta}, \mathcal{L}_{X_{u+\theta}})du + \int_0^{t}\sigma(u,X_{u+\theta}, \mathcal{L}_{X_{u+\theta}})d\Bar{B}_u.
        \end{split}
    \end{equation*}
    Writing in the differential form, we have that
    \begin{equation*}
        \begin{cases}
            dX_{t+\theta} = b(t,X_{t+\theta}, \mathcal{L}_{X_{t+\theta}})dt + \sigma(t,X_{t+\theta}, \mathcal{L}_{X_{t+\theta}})d\Bar{B}_t, \\
            \mathcal{L}_{X_\theta} = \mu.
        \end{cases}
    \end{equation*}
    Then, $\{X_t\}_{t\ge 0}$ and $\{X_{t+\theta}\}_{t\ge 0}$ satisfy the DDSDE (\ref{DDSDE}) with Brownian motions  $\{B_t\}_{t\geq0}$ and  $\{\Bar{B}_t\}_{t\geq0}$, respectively. Therefore, by the weak uniqueness of solutions, we conclude that for any sequence $0 \leq t_1 < t_2< \cdots< t_n$, $(X_{t_1}, X_{t_2},\dots, X_{t_n})$ and $(X_{t_1 + \theta},X_{t_2 + \theta}, \dots, X_{t_n+\theta})$ have the same joint distribution.\hfill $\square$

\noindent {\it Proof of Theorem \ref{thm1}.} By Lemma \ref{lemma}, it suffices to show that there exists a $\theta$-periodic measure  in ${\cal P}_{\vartheta}(\mathbb{R}^d)$ for the DDSDE (\ref{DDSDE}).

First, we assume that condition (\textbf{H3a}) or condition (\textbf{H3b}) holds.     Let $X^{\delta_0}_t$ be the solution to the DDSDE (\ref{DDSDE}) with starting point $0$ and $V$ be a function satisfying condition (\textbf{H3a}) or condition (\textbf{H3b}). By Ito's formula, we get
\begin{eqnarray*}
    V(t,X^{\delta_0}_t)-V(0,0)=\int_0^t(LV)(s,X^{\delta_0}_s,{\cal L}_{X^{\delta_0}_s})ds+\int_0^t\sigma(s,X^{\delta_0}_s,{\cal L}_{X^{\delta_0}_s})\cdot\partial_xV(s,X^{\delta_0}_s)dB_s.
\end{eqnarray*}
Hence, we have that
    \begin{equation}\label{Dec19a}
    -C_0 - V(0,0) \le E[V(t,X^{\delta_0}_t)]-V(0,0)=E\left[\int_0^t(LV)(s,X^{\delta_0}_s,{\cal L}_{X^{\delta_0}_s})ds\right].
 \end{equation}

 (i) If  condition (\textbf{H3a}) holds, then
    $$
    (LV)(s,X^{\delta_0}_s,{\cal L}_{X^{\delta_0}_s})\le -C_1|X^{\delta_0}_s|^{r}+C_2(1+\|{\cal L}_{X^{\delta_0}_s}\|^{\vartheta}_{\vartheta}).
    $$
    Thus, by (\ref{Dec19a}), we deduce that for any $N\in\mathbb{N}$,
    \begin{eqnarray}\label{Oct9bb}
    C_1\int_0^t\|{\cal L}_{X^{\delta_0}_s}\|^{r}_{r}ds&\le& C_2\int_0^t\|{\cal L}_{X^{\delta_0}_s}\|^{\vartheta}_{\vartheta}ds+C_2t+C_0+V(0,0)\nonumber\\
    &\le&\frac{C_2}{N^{r-\vartheta}}\int_0^t\|{\cal L}_{X^{\delta_0}_s}\|^{r}_{r}ds+C_2(1+N^{\vartheta})t+C_0+V(0,0).
    \end{eqnarray}
    We choose an $N_0\in\mathbb{N}$  satisfying $\frac{C_2}{N_0^{r-\vartheta}}\le\frac{C_1}{2}$. Hence, by (\ref{Oct9bb}), we get
    \begin{eqnarray}\label{Oct9b}
    \frac{1}{t}\int_0^t\|{\cal L}_{X^{\delta_0}_s}\|^{r}_{r}ds\le \frac{2\left[C_2(1+N_0^{\vartheta})+\frac{C_0+V(0,0)}{t}\right]}{C_1}.
    \end{eqnarray}

 (ii) If  condition (\textbf{H3b}) holds, then
    $$
    (LV)(s,X^{\delta_0}_s,{\cal L}_{X^{\delta_0}_s})\le -C_1|X^{\delta_0}_s|^{r}+C_2\|{\cal L}_{X^{\delta_0}_s}\|^{r}_{r}+C_3.
    $$
    Thus,  by (\ref{Dec19a}), we get
$$
    C_1\int_0^t\|{\cal L}_{X^{\delta_0}_s}\|^{r}_{r}ds\le C_2\int_0^t\|{\cal L}_{X^{\delta_0}_s}\|^{r}_{r}ds+C_3t+C_0+V(0,0),
$$
which implies that
        \begin{eqnarray}\label{Oct9bbb}
    \frac{1}{t}\int_0^t\|{\cal L}_{X^{\delta_0}_s}\|^{r}_{r}ds\le \frac{C_3+\frac{C_0+V(0,0)}{t}}{C_1-C_2}.
    \end{eqnarray}

By (\ref{Oct9b}) and (\ref{Oct9bbb}), we conclude that if condition (\textbf{H3a}) or condition (\textbf{H3b}) holds then
$$
    \lim_{R\rightarrow\infty}\limsup_{t\rightarrow\infty}\frac{1}{t}\int_0^tP(0,0,s,{\overline{B_{R}}}^c)ds=0,
$$
where $B_R:=\{x\in\mathbb{R}^d: |x|<R\}$ and $P(s,x, t, A)$ is the transition probability function of $\{X_t\}_{t\ge0}$ for $0\leq s<t < \infty$, $x \in \mathbb{R}^d$ and $A\in \mathcal{B}(\mathbb{R}^d)$. Hence, there exist a sequence of natural numbers $T_n\uparrow\infty$ and $\nu\in{\cal P}(\mathbb{R}^d)$ such that
    \begin{equation}\label{Oct9c}
    \nu_n:=\frac{1}{T_n}\int_0^{T_n}P^*_{s}\delta_0ds=\frac{1}{T_n}\int_0^{T_n}P(0,0,s,\cdot)ds\xrightarrow[]{w}\nu\ \ {\rm as}\ n\rightarrow\infty.
    \end{equation}
Moreover, by (\ref{Oct9b}) and (\ref{Oct9bbb}), we get
    $$
    \lim_{R\rightarrow\infty}\limsup_{n\rightarrow\infty}\int_{|y|> R}|y|^{\vartheta}\nu_n(dy)=0.
    $$
    Then, by \cite[Definition 6.8 and Theorem 6.9]{villani2016optimal}, we deduce that $\nu\in {\cal P}_{\vartheta}(\mathbb{R}^d)$ and
    \begin{equation}\label{Oct9d}
    \lim_{n\rightarrow\infty}W_{\vartheta}(\nu_n,\nu)=0.
    \end{equation}

    Let $G$ be an open subset of $\mathbb{R}^d$. We have
    \begin{eqnarray*}
    (P_{\theta}^*\nu_n)(G)&=&\frac{1}{T_n}\int_0^{T_n}(P_{0,\theta}^*P^*_{0,s}\delta_0)(G)ds\\
    &=&\frac{1}{T_n}\int_0^{T_n}(P_{s,s+\theta}^*P^*_{0,s}\delta_0)(G)ds\\
    &=&\frac{1}{T_n}\int_0^{T_n}(P_{s+\theta}^*\delta_0)(G)ds\\
    &=&\frac{1}{T_n}\int_{\theta}^{T_n+\theta}(P_{s}^*\delta_0)(G)ds\\
    &=&\frac{1}{T_n}\int_{\theta}^{T_n+\theta}P(0,0,s,G)ds.
    \end{eqnarray*}
    Then, by (\ref{Oct9c}) and \cite[Theorem 3.3.1]{KurtzMarkov}, we get
    $$
    \liminf_{n\rightarrow\infty}\,(P_{\theta}^*\nu_n)(G)\ge \nu(G).
    $$
    Since $G$ is arbitrary, using \cite[Theorem 3.3.1]{KurtzMarkov} again we get
    \begin{equation}\label{Oct9e}
    P_{\theta}^*\nu_n\xrightarrow[]{w}\nu\ \ {\rm as}\ \ n\rightarrow\infty.
    \end{equation}
   Further, by condition \textbf{(H2)} and (\ref{Oct9d}), we obtain
    $$
    \lim_{n\rightarrow\infty}W_{\vartheta}(P_{\theta}^*\nu_n,P_{\theta}^*\nu)=0,
    $$
    which together with (\ref{Oct9e}) implies that
    $    P_{\theta}^*\nu=\nu$.     We have obtained a $\theta$-periodic measure  in ${\cal P}_{\vartheta}(\mathbb{R}^d)$ for the DDSDE (\ref{DDSDE}).
\vskip 0.2cm
Next, we assume that condition (\textbf{H3c}) holds. We will use the Schauder fixed point theorem to complete the proof. First, let us recall the Schauder fixed point theorem.

\begin{thm}\label{schauder} (cf. \cite[Theorem 11.1.2]{Subrahmanyam}) Let $X$ be a  normed space and $K \subset X$ be a non-empty, compact and convex  subset. Then, every continuous mapping $T: K\rightarrow K$ has a fixed point.
\end{thm}

Let ${\cal M}_1(\mathbb{R}^d)$ be the set of all finite
    signed measures on $\mathbb{R}^d$ with $\int_{\mathbb{R}^d}|x|\mu(dx)<\infty$. It is known that (cf. \cite[page 249 and Corollary 5.4]{CARMONA}) ${\cal M}_1(\mathbb{R}^d)$ is a normed space when equipped with the Kantorovich-Rubinstein norm:
    $$
    \|\mu\|_{{\rm KR}}:=|\mu(\mathbb{R}^d)|+\sup \left\{\int_{\mathbb{R}^d}l(x)\mu(dx): l\in {\rm Lip}_1(\mathbb{R}^d),\,l(0)=0\right\}.
    $$
    Moreover,
\begin{eqnarray}\label{Jan14b}
    \|\mu_1-\mu_2\|_{{\rm KR}}=W_1(\mu_1,\mu_2),\ \ \ \ \mu_1,\mu_2\in {\cal P}_1(\mathbb{R}^d).
\end{eqnarray}

   Let $V$ be a function satisfying condition (\textbf{H3c}).   Define
    $$
    {\cal K}:=\left\{\mu\in{\cal P}(\mathbb{R}^d):\int_{\mathbb{R}^d}V(0,x)\mu(dx)\le |V(0,0)|\right\}.
    $$
Obviously, $\delta_0\in {\cal K}$ and ${\cal K}$ is a convex set. By condition (\textbf{H3c}), we get
\begin{eqnarray}\label{Jan14a}
        \int_{\mathbb{R}^d}|x|^r\mu(dx) \le \frac{1}{C_1}\left\{|V(0,0)|+C_2\right\},\ \ \ \ \forall \mu\in{\cal K}.
\end{eqnarray}
Note that $r>1$. Hence, by (\ref{Jan14b}), (\ref{Jan14a}) and \cite[Definition 6.8 and Theorem 6.9]{villani2016optimal},  we deduce that ${\cal K}$ is a compact subset  of $({\cal M}_1(\mathbb{R}^d),\|\cdot\|_{{\rm KR}})$.

Let $\mu\in {\cal K}$. Suppose that ${\cal L}_{X_0}=\mu$. By Ito's formula, we get
\begin{eqnarray*}
    V(t,X_t)-V(0,X_0)=\int_0^t(LV)(s,X_s,{\cal L}_{X_s})ds
+\int_0^t\sigma(s,X_s,{\cal L}_{X_s})\cdot\partial_xV(s,X_s)dB_s,
\end{eqnarray*}
which together with condition \textbf{(H3c)} implies that
    $$
    E[V(0,X_{\theta})]\le E[V(\theta,X_{\theta})]\le  E[V(0,X_0)].
    $$
Then,
 \begin{eqnarray*}
        \int_{\mathbb{R}^d}V(0,x)P_{\theta}^*\mu(dx)&=& \int_{\mathbb{R}^d}V(0,x){\cal L}_{X_{\theta}}(dx)\\
&\le&\int_{\mathbb{R}^d}V(0,x){\cal L}_{X_0}(dx)\\
& =& \int_{\mathbb{R}^d}V(0,x)\mu(dx)\\
&\le&|V(0,0)|.
\end{eqnarray*}
Thus, $P_{\theta}^*\mu \in {\cal K}$. Since $\mu\in {\cal K}$ is arbitrary, we conclude that $
    P^*_{\theta}{\cal K}\subset {\cal K}$.

Assume that $\mu_n,\mu_0\in {\cal K}$ and $\mu_n$ converges to $\mu_0$ with respect to the $\|\cdot\|_{{\rm KR}}$-norm as $n\rightarrow\infty$. Since $r>\vartheta$, by (\ref{Jan14a}) and  \cite[Definition 6.8 and Theorem 6.9]{villani2016optimal}, we deduce that
    $\lim_{n\rightarrow\infty}W_{\vartheta}(\mu_n,\mu_0)$ $=0$. Then, by condition \textbf{(H2)}, we get $    \lim_{n\rightarrow\infty}W_{\vartheta}(P_{\theta}^*\mu_n,P_{\theta}^*\mu_0)=0$, which implies that $    \lim_{n\rightarrow\infty}$ $W_1(P_{\theta}^*\mu_n,P_{\theta}^*\mu_0)=0$.
    Thus, $P_{\theta}^*\mu_n$ converges to $P_{\theta}^*\mu_0$ with respect to the $\|\cdot\|_{{\rm KR}}$-norm as $n\rightarrow\infty$ by (\ref{Jan14b}). Hence, $P_{\theta}^*$ is a continuous mapping from ${\cal K}$ to itself. Therefore, by Theorem \ref{schauder}, there exists $\Bar{\mu}\in {\cal K}\subset{\mathcal{P}_{\vartheta}(\mathbb{R}^d)}$ such that $P_{\theta}^*\Bar{\mu} = \Bar{\mu}$.{\it Proof.}

\begin{rem} From the proof of Theorem  \ref{thm1}, we can see that under the assumptions of Theorem  \ref{thm1} there exists a $\theta$-periodic measure in ${\cal P}_{\vartheta}(\mathbb{R}^d)$ for the DDSDE (\ref{DDSDE}), and under the assumptions of Corollary  \ref{cor} there exists an invariant measure in ${\cal P}_{\vartheta}(\mathbb{R}^d)$ for the DDSDE (\ref{timeIndDDSDE}).
\end{rem}

\section{Examples}\label{sec3}\setcounter{equation}{0}

In this section, we provide applications for the criteria presented in Section 1. We will mainly use Theorem \ref{thm1} to give examples for the existence of periodic solutions to DDSDEs. However,  by assuming that the coefficients $b$ and $\sigma$ are time-independent and modifying the assumptions correspondingly, we can use Corollary \ref{cor} to give similar examples for the existence of stationary solutions to DDSDEs.

Let $\vartheta \in [1,\infty)$.  We put the following assumptions:
\begin{itemize}
    \item[(\textbf{A1})] $\sigma(t,0,\delta_0)$ is bounded on $[0,\theta)$ and  there exist $K_{\sigma,1},K_{\sigma,2}>0$ such that
\begin{eqnarray*}
&&\|\sigma(t,x,\mu)-\sigma(t,y,\nu)\|^2\le K_{\sigma,1}|x-y|^2+K_{\sigma,2}\{W_\vartheta(\mu,\nu)\}^2,\\
&&\ \ \ \ \ \ \ \  \forall t\in[0,\theta), x,y\in\mathbb{R}^d, \mu,\nu\in{\cal P}_{\vartheta}(\mathbb{R}^d).
\end{eqnarray*}
Additionally, $K_{\sigma,2}=0$ when $\vartheta\in[1,2)$.
    \item[(\textbf{A2})]  $b$ is locally bounded on $[0,\infty)\times\mathbb{R}^d\times {\cal P}_{\vartheta}(\mathbb{R}^d)$ and  there exist $K_{b,1},K_{b,2},K_{b,3}>0$ such that
\begin{eqnarray*}
&&\langle b(t,x,\mu)-b(t,y,\nu),x-y\rangle\le K_{b,1}|x-y|^2+K_{b,2}W_\vartheta(\mu,\nu)|x-y|,\\
&&\ \ \ \ \ \ \ \ \forall t\in[0,\theta), x,y\in\mathbb{R}^d, \mu,\nu\in{\cal P}_{\vartheta}(\mathbb{R}^d),
\end{eqnarray*}
and
\begin{eqnarray}\label{Dec25a}
|b(t,0,\mu)|\le K_{b,3}(1+\|\mu\|_{\vartheta}),\quad\forall t\in[0,\theta), \mu\in{\cal P}_{\vartheta}(\mathbb{R}^d).
\end{eqnarray}
\end{itemize}

\begin{rem}\label{rem311}
Note that conditions (H0), (A1) and (A2) imply that the coefficients $b$ and $\sigma$ satisfy the continuity, monotonicity and growth conditions of \cite[Theorem 2.1]{wang2017distributiondependent}. Then, by \cite[Theorem 2.1]{wang2017distributiondependent}, we conclude that conditions (H1), (H2) hold and $\int_0^t \|{\cal L}_{X^{\delta_0}_s}\|^{\vartheta}_{\vartheta}ds<\infty$ for any  $t>0$. Moreover, since $\|\mu\|_{\vartheta}\le \|\mu\|_{\vartheta'}$ for any $\mu\in \prd$ and $\vartheta'>\vartheta$,  conditions (H0), (A1) and (A2) also imply that condition (H2) holds with $\vartheta$ replaced by any $\vartheta'>\vartheta$ and $\int_0^t \|{\cal L}_{X_s}\|^{\vartheta'}_{\vartheta'}ds<\infty$ for any ${\cal L}_{X_0}\in {\cal P}_{\vartheta'}(\mathbb{R}^d)$ and $t>0$.
\end{rem}
\vskip 0.2cm
\begin{exa}\label{example 3.1}\end{exa} First, we present a general result.

\vskip 0.2cm
\begin{pro}\label{proDec12}
Suppose that  the coefficients $b$ and $\sigma$ satisfy conditions \textbf{(H0)}, \textbf{(A1)}, \textbf{(A2)} and one of the following two conditions:
\begin{itemize}
 \item[(\textbf{A3})]  $\vartheta\in[1,2)$ and there exist $K_{b,4}>\frac{K_{\sigma,1}}{2}$, $K_{b,5}>0$  such that
$$
\langle b(t,x,\mu),x\rangle\le -K_{b,4}|x|^2+K_{b,5}(1+\|\mu\|_{\vartheta}^{\vartheta}),\quad
\forall t\in[0,\theta), x\in\mathbb{R}^d, \mu\in{\cal P}_{\vartheta}(\mathbb{R}^d).
$$
    \item[(\textbf{A4})]  $\vartheta\in[2,\infty)$ and there exist $K_{b,4},K_{b,5}>0$, $r>\vartheta$  such that
$$
\langle b(t,x,\mu),x\rangle\le -K_{b,4}|x|^r+K_{b,5}(1+\|\mu\|_{\vartheta}^{\vartheta}),\quad
\forall t\in[0,\theta), x\in\mathbb{R}^d, \mu\in{\cal P}_{\vartheta}(\mathbb{R}^d).
$$
\end{itemize}
Then, there exists a $\theta$-periodic solution to the DDSDE (\ref{DDSDE}).
\end{pro}

\noindent {\it Proof.} By Theorem \ref{thm1} and  Remark \ref{rem311}, to show that there exists a $\theta$-periodic solution to the DDSDE (\ref{DDSDE}), we need only show that there exists  $V\in C^{1,2}([0,\infty)\times\mathbb{R}^d)$ satisfying condition (H3a).

Let $V(t,x)=|x|^2$ for $t\ge0$ and $x\in\mathbb{R}^d$. By (\ref{Lya Generator}), we get
 \begin{eqnarray}\label{Dec11a}
                    (LV)(t,x,\mu)=2\langle b(t,x,\mu), x\rangle+{\rm trace}\left(\sigma\sigma^T(t,x,\mu)\right).
                \end{eqnarray}
First, we assume that condition (\textbf{A3}) holds.
Note that condition (A1) implies that
\begin{eqnarray*}
&&\|\sigma(t,x,\mu)\|\le K^{\frac{1}{2}}_{\sigma,1}|x|+\sup_{0\le s<\theta}\|\sigma(s,0,\delta_0)\|,\quad\forall t\in[0,\theta), x\in\mathbb{R}^d, \mu\in{\cal P}_{\vartheta}(\mathbb{R}^d).
\end{eqnarray*}
Then, by (\ref{Dec11a}) and condition (A3), we get
\begin{eqnarray*}
(LV)(t,x,\mu)&\le& -2K_{b,4}|x|^{2}+2K_{b,5}(1+\|\mu\|_{\vartheta}^{\vartheta})+\left(K^{\frac{1}{2}}_{\sigma,1}|x|+\sup_{0\le s<\theta}\|\sigma(s,0,\delta_0)\|\right)^2\\
&\le&-2K_{b,4}|x|^{2}+2K_{b,5}(1+\|\mu\|_{\vartheta}^{\vartheta})+K_{\sigma,1}|x|^2+\left(\sup_{0\le s<\theta}\|\sigma(s,0,\delta_0)\|\right)^2\\
&&+\frac{2K_{b,4}-K_{\sigma,1}}{2}|x|^2+\frac{2K_{\sigma,1}(\sup_{0\le s<\theta}\|\sigma(s,0,\delta_0)\|)^2}{2K_{b,4}-K_{\sigma,1}}\\
&=&-\frac{2K_{b,4}-K_{\sigma,1}}{2}|x|^2+2K_{b,5}(1+\|\mu\|_{\vartheta}^{\vartheta})\\
&&+\frac{(2K_{b,4}+K_{\sigma,1})(\sup_{0\le s<\theta}\|\sigma(s,0,\delta_0)\|)^2}{2K_{b,4}-K_{\sigma,1}}.
\end{eqnarray*}
Thus, condition (H3a) holds with $r=2$.
\vskip 0.2cm
Next, we assume that condition (\textbf{A4}) holds. Note that in this case condition (A1) implies that
\begin{eqnarray*}
&&\|\sigma(t,x,\mu)\|\le \left(K_{\sigma,1}|x|^2+K_{\sigma,2}\|\mu\|_{\vartheta}^2\right)^{\frac{1}{2}}+\sup_{0\le s<\theta}\|\sigma(s,0,\delta_0)\|,\nonumber\\
&&\ \ \ \ \ \ \ \  \forall t\in[0,\theta), x\in\mathbb{R}^d, \mu\in{\cal P}_{\vartheta}(\mathbb{R}^d).
\end{eqnarray*}
Then, by (\ref{Dec11a}) and condition (A4), we get
\begin{eqnarray*}
(LV)(t,x,\mu)&\le& -2K_{b,4}|x|^{r}+2K_{b,5}(1+\|\mu\|_{\vartheta}^{\vartheta})\\
&&+\left\{\left(K_{\sigma,1}|x|^2+K_{\sigma,2}\|\mu\|^2_{\vartheta}\right)^{\frac{1}{2}}+\sup_{0\le s<\theta}\|\sigma(s,0,\delta_0)\|\right\}^2\\
&\le&2\bigg\{-K_{b,4}|x|^{r}+K_{\sigma,1}|x|^2+(K_{b,5}+K_{\sigma,2})(1+\|\mu\|^{\vartheta}_{\vartheta})\\
&&\ \ \ \ +\left(\sup_{0\le s<\theta}\|\sigma(s,0,\delta_0)\|\right)^2\bigg\}.
\end{eqnarray*}
Thus, condition (H3a) holds. Therefore, the proof is complete.\hfill $\square$

We now give a concrete example. For $x=(x_1,\dots,x_d)\in\mathbb{R}^d$ and $m\in\mathbb{N}$, denote $x^m=(x_1^m,\dots,x_d^m)$. Let $\gamma_1,\gamma_2,\gamma_3$  be three bounded $\theta$-periodic functions on $[0,\infty)$ with $\inf_{0\le s<\theta}\gamma_1(s)>0$,
\[
    b_0: \mathbb{R}^d\rightarrow \mathbb{R}^d,\quad  \quad \sigma_0: \mathbb{R}^d \rightarrow \mathbb{R}^d\otimes\mathbb{R}^d
\]
be two Lipschitz continuous maps, $n\in\mathbb{N}$, and $\alpha,\beta\in\mathbb{R}$.
For $t\ge0, x\in\mathbb{R}^d, \mu\in\prd$, define
$$
b(t,x,\mu):=-\gamma_1(t)x^{2n+1}+\gamma_2(t)\int_{\mathbb{R}^d}b_0(x-\alpha z)\mu(dz),
$$
and
$$
\sigma(t,x,\mu):=\gamma_3(t)\int_{\mathbb{R}^d}\sigma_0(x-\beta z)\mu(dz).
$$

Set $\vartheta=2$. It is easy to see that condition (H0) and (\ref{Dec25a}) hold. Define
$$
B_0:=|\nabla b_0|_{\infty},\quad C_0:=\sup_{|v|=1,x\in\mathbb{R}^d}|(\nabla_v\sigma_0)(x)|^2,\quad K_0:=\sup_{x\not=y\in\mathbb{R}^d}\frac{\langle b_0(x)-b_0(y),x-y\rangle}{|x-y|^2},
$$
where  $\nabla_v$ denotes the directional derivative along $v$. Then, $B_0,C_0,|K_0|<\infty$ since $b_0$ and $\sigma_0$ are Lipschitz continuous. For any $\pi\in \mathcal{C}(\mu,\nu)$, we have that
\begin{eqnarray*}
&&\langle b(t,x,\mu)-b(t,y,\nu),x-y\rangle\\
&=&-\gamma_1(t)\langle x^{2n+1}-y^{2n+1}, x-y\rangle+\gamma_2(t)\int_{\mathbb{R}^d}\left\{\langle b_0(x-\alpha z)-b_0(y-\alpha z),x-y\rangle\right.\\
&&\ \ \ \ \ \ \ \ \ \ \ \ \ \ \ \ \ \ \ \ \ \left.+\langle b_0(y-\alpha z)-b_0(y-\alpha z'),x-y\rangle\right\}\pi(dz,dz'),
\end{eqnarray*}
which implies that
$$
\langle b(t,x,\mu)-b(t,y,\nu),x-y\rangle\le \sup_{0\le s<\theta}|\gamma_2(s)|\cdot\left(|K_0|\cdot |x-y|^2+B_0|\alpha|W_1(\mu,\nu)|x-y|\right).
$$
Similarly, we get
\begin{eqnarray*}
\|\sigma(t,x,\mu)-\sigma(t,y,\nu)\|^2
&\le& C_0\sup_{0\le s<\theta}|\gamma_3(s)|^2\cdot \left( |x-y|+|\beta|W_1(\mu,\nu)\right)^2\\
&\le& C_0(1+|\beta|)\sup_{0\le s<\theta}|\gamma_3(s)|^2\cdot \left( |x-y|^2+|\beta|\{W_2(\mu,\nu)\}^2\right).
\end{eqnarray*}
Thus, conditions (A1) and (A2) hold. Finally, we have
\begin{eqnarray*}
\langle b(t,x,\mu),x\rangle
&\le&-\inf_{0\le s<\theta}\gamma_1(s)\cdot\sum_{i=1}^dx_i^{2n+2}\\
&&+\sup_{0\le s<\theta}|\gamma_2(s)|\cdot|x|\left\{|b_0(0)|+B_0(|x|+|\alpha|\cdot\|\mu\|_1)\right\}\\
&\le&-\inf_{0\le s<\theta}\gamma_1(s)\cdot\sum_{i=1}^dx_i^{2n+2}\\
&&+\sup_{0\le s<\theta}|\gamma_2(s)|\cdot\left\{ \left[|b_0(0)|+B_0|x|\right]|x|+\frac{B_0|\alpha|}{2}(|x|^2+\|\mu\|_2^2)\right\},
\end{eqnarray*}
which implies that condition (A4) holds with $r=2n+2$. Hence all conditions of Proposition \ref{proDec12} hold and therefore there exists a $\theta$-periodic solution to the DDSDE (\ref{DDSDE}).

\begin{rem}
In the above example, if $\gamma_1>0,\gamma_2,\gamma_3$ are constants, then we can show  that there exists a stationary solution to the DDSDE (\ref{timeIndDDSDE})  by Corollary \ref{cor}. In \cite[Example 5.1]{LIU2022138.MarkovSwitching}, by assuming that $b_0$ and $\sigma_0$ are linear functions, $d=n=1$ and $\alpha=\beta\in(-1,1)$,  Liu and Ma used the Banach fixed point theorem to obtain  the existence and uniqueness of invariant measures for  the DDSDE (\ref{timeIndDDSDE}).
\end{rem}

\begin{exa}\label{exam3.2}\end{exa}  In this example, we consider the time-dependent homogeneous Landau equation with Maxwell molecules. Let $\alpha,\beta\in\mathbb{R}$ and $\gamma_2\ge0,\gamma_3$  be two  bounded $\theta$-periodic functions on $[0,\infty)$ such that
\begin{eqnarray}\label{Dec26a}
\inf_{0\le s<\theta}\gamma_2(s)>|\alpha|\sup_{0\le s<\theta}\gamma_2(s)+\frac{3(1+|\beta|)^2}{2}\sup_{0\le s<\theta}|\gamma_3(s)|^2.
\end{eqnarray}
Define
\[
    b_0(x):=-2x,\quad  \quad \sigma_0(x):= \begin{pmatrix}
    x_2 & 0 & x_3 \\
    -x_1 & x_3 & 0 \\
   0 & -x_2 & -x_1
\end{pmatrix},\quad  \quad x=(x_1,x_2,x_3)\in\mathbb{R}^3,
\]
$$
b(t,x,\mu):=\gamma_2(t)\int_{\mathbb{R}^3}b_0(x-\alpha z)\mu(dz),\quad  \quad t\ge0, x\in\mathbb{R}^3, \mu\in{\cal P}(\mathbb{R}^3),
$$
and
$$
\sigma(t,x,\mu):=\gamma_3(t)\int_{\mathbb{R}^d}\sigma_0(x-\beta z)\mu(dz),\quad  \quad t\ge0, x\in\mathbb{R}^3, \mu\in{\cal P}(\mathbb{R}^3).
$$

Set $\vartheta=2$. Following the argument of Example \ref{example 3.1} with $B_0 = C_0 = 2$ and $K_0 = -2$, we can show that conditions (A1) and (A2) hold. Then, by Theorem \ref{thm1} and  Remark \ref{rem311}, to show that there exists a $\theta$-periodic solution to the DDSDE (\ref{DDSDE}), we need only show that there exists  $V\in C^{1,2}([0,\infty)\times\mathbb{R}^d)$ satisfying condition (H3b).

By (\ref{Dec26a}), we can choose a $q > 1$ such that
\begin{eqnarray}\label{Dec26aaaa}
\inf_{0\le s<\theta}\gamma_2(s)>|\alpha|\sup_{0\le s<\theta}\gamma_2(s)+\frac{(2q+1)(1+|\beta|)^2}{2}\sup_{0\le s<\theta}|\gamma_3(s)|^2.
\end{eqnarray}
Define
       $$
    V(t,x):=|x|^{2q},\ \ \ \ t\ge0,x\in\mathbb{R}^3.
    $$
For $1\le i\le 3$, we have
    $$
    	\partial_{x_i} V(t,x)=2q|x|^{2(q-1)}x_i,
    $$
    and for $1\le i,j\le 3$, we have
\begin{eqnarray}\label{Dec12dd}
\partial^2_{x_ix_j} V(x)=2q|x|^{2(q-2)}\left[2(q-1)x_ix_j+|x|^2\delta_{ij}\right].
\end{eqnarray}
Then,
\begin{eqnarray}\label{Dec12z11}
&&\langle b(t,x,\mu),\partial_x V(t,x)\rangle\nonumber\\
&=& 2q|x|^{2(q-1)}\langle b(t,x,\mu),x\rangle\nonumber\\
&\le&4q|x|^{2q-1}\left(-\inf_{0\le s<\theta}\gamma_2(s)\cdot |x|+|\alpha|\sup_{0\le s<\theta}\gamma_2(s)\cdot\|\mu\|_1\right).
\end{eqnarray}
By (\ref{Dec12dd}) and the Cauchy-Schwarz inequality for the Frobenius inner product, we get
    \begin{eqnarray}\label{Dec12z22}
    &&{\rm trace}\left(\sigma\sigma^T(t,x,\mu)\cdot\partial_x^2V(t,x)\right)\nonumber\\
&\le&{\rm trace}\left(\sigma\sigma^T(t,x,\mu)\right)\cdot \left|{\rm trace}\left(\partial_x^2V(t,x)\right)\right|\nonumber\\
&=&2q(2q+1)|x|^{2(q-1)}\cdot \|\sigma(t,x,\mu)\|^2\nonumber\\
&\le&4q(2q+1)|x|^{2(q-1)}\sup_{0\le s<\theta}|\gamma_3(s)|^2\cdot \left( |x|+|\beta|\cdot\|\mu\|_1\right)^2\nonumber\\
&\le&4q(2q+1)|x|^{2(q-1)}\sup_{0\le s<\theta}|\gamma_3(s)|^2\cdot \left\{(1+|\beta|) |x|^2+(|\beta|+\beta^2)\|\mu\|^2_1\right\}.
    \end{eqnarray}
 Thus,   by (\ref{Lya Generator}),  (\ref{Dec12z11}), (\ref{Dec12z22}) and  Young's inequality for products, we deduce that for $t\ge0, x\in\mathbb{R}^3, \mu\in{\cal P}_{2q}(\mathbb{R}^3)$,
\begin{eqnarray*}
    &&LV(t,x,\mu)\\
    &\le&4q|x|^{2q-1}\left(-\inf_{0\le s<\theta}\gamma_2(s)\cdot |x|+|\alpha|\sup_{0\le s<\theta}\gamma_2(s)\cdot\|\mu\|_1\right)\\
    &&+2q(2q+1)|x|^{2(q-1)}\sup_{0\le s<\theta}|\gamma_3(s)|^2\cdot \left\{(1+|\beta|) |x|^2+(|\beta|+\beta^2)\|\mu\|^2_1\right\}\\
    &\le&-2q\left\{2\inf_{0\le s<\theta}\gamma_2(s)-(2q+1)(1+|\beta|)\sup_{0\le s<\theta}|\gamma_3(s)|^2\right\}|x|^{2q}\\
    &&+4q|\alpha|\sup_{0\le s<\theta}\gamma_2(s)\cdot\|\mu\|_{2q}\cdot |x|^{2q-1}\\
&&+2q(2q+1)(|\beta|+\beta^2)\sup_{0\le s<\theta}|\gamma_3(s)|^2\cdot \|\mu\|^{2}_{2q}\cdot|x|^{2q-2}\\
    &\le&-2q\left\{2\inf_{0\le s<\theta}\gamma_2(s)-(2q+1)(1+|\beta|)\sup_{0\le s<\theta}|\gamma_3(s)|^2\right\}|x|^{2q}\\
    &&+2|\alpha|\sup_{0\le s<\theta}\gamma_2(s)\cdot \|\mu\|^{2q}_{2q}+2(2q-1)|\alpha|\sup_{0\le s<\theta}\gamma_2(s)\cdot |x|^{2q}\\
    &&+2(2q+1)(|\beta|+\beta^2)\sup_{0\le s<\theta}|\gamma_3(s)|^2\cdot  \|\mu\|^{2q}_{2q}\\
&&+2(q-1)(2q+1)(|\beta|+\beta^2)\sup_{0\le s<\theta}|\gamma_3(s)|^2\cdot  |x|^{2q}\\
    &=&-2\Bigg\{2q\inf_{0\le s<\theta}\gamma_2(s)-(2q-1)|\alpha|\sup_{0\le s<\theta}\gamma_2(s)\\
&&\ \ \ \ \ \ \ -(2q+1)(1+|\beta|)[q+(q-1)|\beta|]\sup_{0\le s<\theta}|\gamma_3(s)|^2\Bigg\}|x|^{2q}\\
    &&+2\Bigg\{|\alpha|\sup_{0\le s<\theta}\gamma_2(s)+(2q+1)(|\beta|+\beta^2)\sup_{0\le s<\theta}|\gamma_3(s)|^2\Bigg\} \|\mu\|^{2q}_{2q}.
\end{eqnarray*}
Hence, condition (H3b) holds with $r=2q$ by (\ref{Dec26aaaa}). Therefore, there exists a $\theta$-periodic solution to the DDSDE (\ref{DDSDE}) by Theorem \ref{thm1}.

\begin{rem}
In Example \ref{exam3.2}, if $|\alpha|<1$, $\beta\in\mathbb{R}$, and $\gamma_2$ and $\gamma_3$ are constants such that
$$
\gamma_2>\frac{3(1+|\beta|)^2}{2(1-|\alpha|)}\gamma_3^2,
$$
then we can show that there exists a stationary solution to the DDSDE (\ref{timeIndDDSDE}) by Corollary \ref{cor}. In \cite[Corollary 6.2]{wang2017distributiondependent}, by assuming that $\gamma_2=\gamma_3=1$ and $2(|\alpha| + |\beta|) + \beta^2 < 1$,  Wang used the Banach fixed point theorem to obtain the existence and uniqueness of invariant measures for the DDSDE (\ref{timeIndDDSDE}).
\end{rem}

\begin{exa}\end{exa} In this example, we explain how to use condition (H3c) to establish the existence of periodic solutions to DDSDEs.

\noindent (a)  Suppose that  the coefficients $b$ and $\sigma$ satisfy conditions \textbf{(H0)}, \textbf{(A1)}, \textbf{(A2)} and the following condition:
\begin{itemize}
    \item[(\textbf{A5})]  There exists $q>\max\{1,\frac{\vartheta}{2}\}$ such that
\begin{eqnarray}\label{Dec12f}
\langle b(t,x,\mu),x\rangle+\left(q-1+\frac{d}{2}\right)\|\sigma(t,x,\mu)\|^2\le0,\quad \forall t\ge0,  x\in\mathbb{R}^d, \mu\in{\cal P}(\mathbb{R}^d).
\end{eqnarray}
\end{itemize}
Then, there exists a $\theta$-periodic solution to the DDSDE (\ref{DDSDE}).

 In fact, by Theorem \ref{thm1}  and  Remark \ref{rem311}, we need only show that condition (H3c) holds. Let $q>\max\{1,\frac{\vartheta}{2}\}$ such that condition (A5) holds. Define
    $$
    g(y):=y-e^{-y},\ \ \ \ y\ge0,
    $$
and
    $$
    V(t,x):=g(|x|^{2q}),\ \ \ \ t\ge0,x\in\mathbb{R}^d.
    $$

It is easy to see that the first part of condition (H3c) holds. For $1\le i\le d$, we have
    $$
    	\partial_{x_i} V(t,x)=2q\left(1+e^{-|x|^{2q}}\right)
    		|x|^{2(q-1)}x_i,
    $$
    and for $1\le i,j\le d$, we have
\begin{eqnarray}\label{Dec12d}
\partial^2_{x_ix_j} V(x)&=&-4q^2e^{-|x|^{2q}}
    		|x|^{4(q-1)}x_ix_j\nonumber\\
&&+2q\left(1+e^{-|x|^{2q}}\right)|x|^{2(q-2)}[2(q-1)x_ix_j+|x|^2\delta_{ij}].
\end{eqnarray}
Then,
\begin{eqnarray}\label{Dec12z1}
\langle b(t,x,\mu),\partial_x V(t,x)\rangle= 2q\left(1+e^{-|x|^{2q}}\right)
    		|x|^{2(q-1)}\langle b(t,x,\mu),x\rangle.
\end{eqnarray}
By (\ref{Dec12d}) and the Cauchy-Schwarz inequality for the Frobenius inner product, we get
    \begin{eqnarray}\label{Dec12z2}
    &\ \ \ \ \ \ \ \ &{\rm trace}\left(\sigma\sigma^T(t,x,\mu)\cdot\partial_x^2V(t,x)\right)\nonumber\\
&\ \ \ \ \ \ \ \ \le&{\rm trace}\left(\sigma\sigma^T(t,x,\mu)\right)\cdot \left|{\rm trace}\left(\partial_x^2V(t,x)\right)\right|\nonumber\\
&\ \ \ \ \ \ \ \ =&\|\sigma(t,x,\mu)\|^2\cdot \left|-4q^2e^{-|x|^{2q}}
    		|x|^{2(2q-1)}+2q[2(q-1)+d]\left(1+e^{-|x|^{2q}}\right)|x|^{2(q-1)}\right|\nonumber\\
&\ \ \ \ \ \ \ \ \le& 2q[2(q-1)+d]\left(1+e^{-|x|^{2q}}\right)|x|^{2(q-1)}\|\sigma(t,x,\mu)\|^2.
    \end{eqnarray}
 Thus,   by (\ref{Lya Generator}), (\ref{Dec12f}), (\ref{Dec12z1}) and (\ref{Dec12z2}), we deduce that
    \begin{eqnarray*}
    LV(t,x,\mu)&\le&2q\left(1+e^{-|x|^{2q}}\right)
    		|x|^{2(q-1)}\left\{\langle b(t,x,\mu),x\rangle+\left(q-1+\frac{d}{2}\right)\|\sigma(t,x,\mu)\|^2\right\}\\
&\le & 0,\ \ \ \ \forall t\ge0, x\in\mathbb{R}^d, \mu\in{\cal P}(\mathbb{R}^d).
    \end{eqnarray*}
 Therefore,  the second part of condition (H3c) holds.

\begin{rem} Some other Lyapunov
conditions of type (\textbf{A5}) have been used to establish the existence of invariant measures for DDSDEs. We refer the reader to \cite[condition (H$_3$)]{Bao2022} and \cite[condition (H1)]{Zhang.nonuniqueness}.
\end{rem}

\noindent (b)  Suppose that  the coefficients $b$ and $\sigma$ satisfy conditions \textbf{(H0)}, \textbf{(A1)}, \textbf{(A2)} and the following condition:
\begin{itemize}
    \item[(\textbf{A6})]  $\langle b(t,x,\mu),x\rangle \le 0$ for any $t\ge0,  x\in\mathbb{R}^d, \mu\in{\cal P}(\mathbb{R}^d)$ and there exists $C_{\sigma}>0$ such that $\sigma(t,x,\mu)=0$ for  any $t\ge0,  |x|>C_{\sigma}, \mu\in{\cal P}(\mathbb{R}^d)$.
\end{itemize}
Then, there exists a $\theta$-periodic solution to the DDSDE (\ref{DDSDE}).

     In fact, by Theorem \ref{thm1} and  Remark \ref{rem311}, we need only show that condition (H3c) holds. Let $q>\max\{1,\frac{\vartheta}{8}\}$. Define
    \begin{eqnarray*}
    	g(y):=\begin{cases}
    		0, &0\le y\le C_{\sigma}^{2q},\\
    		(y-C_\sigma^{2q})^4, &		y> C_{\sigma}^{2q},
    	\end{cases}
    	\end{eqnarray*}
    and
    $$
    V(t,x):=g(|x|^{2q}),\ \ \ \ t\ge0,x\in\mathbb{R}^d.
    $$

It is easy to see that the first part of condition (H3c) holds. For $1\le i\le d$, we have
    \begin{eqnarray*}
    	\partial_{x_i} V(t,x)=\begin{cases}
    		0, &|x|\le C_{\sigma},\\
    		8q(|x|^{2q}-C_{\sigma}^{2q})^3|x|^{2(q-1)}x_i, &		|x|> C_{\sigma},
    	\end{cases}
    	\end{eqnarray*}
    and for $1\le i,j\le d$, we have
    \begin{eqnarray*}
    \partial^2_{x_ix_j} V(t,x)=\begin{cases}
    		0, &|x|\le C_{\sigma},\\
    		8q(|x|^{2q}-C_{\sigma}^{2q})^2|x|^{2(q-2)}\\
            \times\left\{(|x|^{2q}-C_{\sigma}^{2q})[|x|^2\delta_{ij}+2(q-1)x_ix_j]+6q|x|^{2q}x_ix_j\right\}, &		|x|> C_{\sigma}.
    \end{cases}
    \end{eqnarray*}
Then, by (\ref{Lya Generator}) and condition (A6), we get
    $$
    (LV)(t,x,\mu)\le 0,\ \ \ \ \forall t\ge0, x\in\mathbb{R}^d, \mu\in{\cal P}(\mathbb{R}^d).
    $$
 Therefore,  the second part of condition (H3c) holds.

\begin{rem} Note that if $b(t,x,\mu)=-\frac{x}{|x|^2}$ when $|x|$ is large, condition \textbf{(A3)} or condition \textbf{(A4)} is not satisfied. However, if condition \textbf{(A6)} is satisfied, there exists a $\theta$-periodic solution to the DDSDE (\ref{DDSDE}).
\end{rem}

\end{document}